\documentclass{amsart}

\usepackage{xypic}
\input xy
\xyoption{all}
\usepackage{epsfig}
\usepackage{amsthm}
\usepackage{amssymb}
\usepackage{amsmath}
\usepackage{amscd}
\usepackage{color}

\newcommand{\lb}{\varLambda}

\newcommand{\A}{\mathcal{A}}
\def \<{\langle}
\def \>{\rangle}

\newcommand{\bg}{\begin{equation}}
\newcommand{\ed}{\end{equation}}
\newcommand{\bga}{\begin{eqnarray}}
\newcommand{\eda}{\end{eqnarray}}

\def\cbdu{\par{\raggedleft$\Box$\par}}

\newtheorem {Theorem}  {Theorem}

\numberwithin{Theorem}{section}

\newtheorem {Lemma}[Theorem]  {Lemma}

\theoremstyle{definition}

\theoremstyle{remark}

\def \l{\lambda}
\expandafter\chardef\csname pre amssym.def
at\endcsname=\the\catcode`\@ \catcode`\@=11
\def\undefine#1{\let#1\undefined}
\def\newsymbol#1#2#3#4#5{\let\next@\relax
 \ifnum#2=\@ne\let\next@\msafam@\else
 \ifnum#2=\tw@\let\next@\msbfam@\fi\fi
 \mathchardef#1="#3\next@#4#5}
\def\mathhexbox@#1#2#3{\relax
 \ifmmode\mathpalette{}{\m@th\mathchar"#1#2#3}%
 \else\leavevmode\hbox{$\m@th\mathchar"#1#2#3$}\fi}
\def\hexnumber@#1{\ifcase#1 0\or 1\or 2\or 3\or 4\or 5\or 6\or 7\or 8\or
 9\or A\or B\or C\or D\or E\or F\fi}

\font\teneufm=eufm10 \font\seveneufm=eufm7 \font\fiveeufm=eufm5
\newfam\eufmfam
\textfont\eufmfam=\teneufm \scriptfont\eufmfam=\seveneufm
\scriptscriptfont\eufmfam=\fiveeufm

\catcode`\@=\csname pre amssym.def at\endcsname

\newcounter{remark}
\def \grad {\nabla}
\def \ls {\lambda_q^{2s}}

\newcommand{\el}{\varepsilon}

\renewcommand{\k}{\kappa}
\renewcommand{\l}{\lambda}

\newcommand{\R}{\mathbf{R}}

\renewcommand{\div}{\mbox{div}}

\def  \R   {{\mathbb R}}
\def  \Z   {{\mathbb Z}}

\def  \T   {{\mathbb T}}

\def  \12  {{\frac{1}{2}}}
\def  \p   {\partial}

\def\build#1_#2^#3{\mathrel{\mathop{\kern 0pt#1}\limits_{#2}^{#3}}}

\begin{document}

\title[Determining Wavenumber for 3D NSE]{Kolmogorov's dissipation number and the number of degrees of freedom for the 3D Navier-Stokes equations}

\author [Alexey Cheskidov]{Alexey Cheskidov}
\address{Department of Mathematics, Stat. and Comp. Sci.,  University of Illinois Chicago, Chicago, IL 60607,USA}
\email{acheskid@uic.edu} 
\author [Mimi Dai]{Mimi Dai}
\address{Department of Mathematics, Stat. and Comp. Sci.,  University of Illinois Chicago, Chicago, IL 60607,USA}
\email{mdai@uic.edu}

\thanks{A. Cheskidov was partially supported by the NSF Grant
DMS--1517583 and M. Dai was partially supported by the NSF Grant DMS--1815069.}





\begin{abstract}
Kolmogorov's theory of turbulence predicts that only wavenumbers below some critical value, called Kolmogorov's dissipation
number, are essential to describe the evolution of a three-dimensional fluid flow. A determining wavenumber, first introduced by
Foias and Prodi for the 2D Navier-Stokes equations, is a mathematical analog of Kolmogorov's number. 
The purpose of this paper is to prove the existence of a time-dependent determining wavenumber for the 3D Navier-Stokes equations whose time
average is bounded by Kolmogorov's dissipation wavenumber for all solutions on the global attractor whose  intermittency
is not extreme.

\bigskip


\hspace{0.02cm}CLASSIFICATION CODE: 35Q35, 37L30.\\

\end{abstract}

\maketitle

\section{Introduction}

The Navier-Stokes equations (NSE) on a torus $\mathbb{T}^3$  are given by 
\begin{equation}
  \label{nse}
  \left\{
    \begin{array}{l}
      u_t + (u \cdot \grad) u - \nu\Delta u + \grad p = f \\
      \grad \cdot u = 0,
    \end{array}
  \right.
\end{equation}
where $u$ is the velocity, $p$ is the pressure, and $f$ is the external force. We assume that
$f$ has zero mean, and consider zero mean solutions. We also assume that $f\in H^{-1}$ or $f$ is translationally bounded in $L^2_{\mathrm{loc}}(\mathbb{R},H^{-1})$.

In this paper we investigate the number of degrees of freedom of a three-dimensional fluid flow governed by \eqref{nse}. Kolmogorov's theory of turbulence \cite{K41}
predicts that there is a wavenumber $\kappa_{\mathrm{d}}$ above which the viscous forces dominate. This suggests that the frequencies above $\kappa_{\mathrm{d}}$
should not affect the dynamics and the number of degrees of freedom is of order $\kappa_{\mathrm{d}}^3$. A natural question is whether this can be justified mathematically.

The notion of determining modes, which allows us to define the degrees of freedom mathematically, was introduced by  Foias and Prodi in \cite{FP} where they showed that high modes of
a solution to the 2D NSE are controlled by low modes asymptotically as time goes to infinity.  Then the number of these
determining modes was estimated by Foias, Manley, Temam, and Treve \cite{FMTT} and later improved by Jones and Titi  \cite{JT}.
We refer the readers to  \cite{CFMT,CFT,FJKT,FJKTgeneral,FMRT,FT,FT84, FT-attractor,FTiti} and references therein for more background and related results.

In this paper we are concerned with 3D flows governed by \eqref{nse}, for which the existence of regular solutions is one of the Millennium open questions. Therefore,
we study weak solutions, whose existence was proved by Leray \cite{L34}. In \cite{CDK}, Cheskidov, Dai, and Kavlie  proved the existence of a determining wavenumber $\lb_u(t)$, defined for each individual trajectory
$u(t)$, whose average is uniformly bounded on the global attractor. More precisely, it was shown  that two solutions $u(t)$ and $v(t)$ on the global attractor are identical, provided
their projections below modes $\max\{\lb_u, \lb_v\}$ coincide. This recovered the results by Constantin, Foias, Manley, and Temam \cite{CFMT} in the 
case where $\| \nabla u(t)\|_{L^2}^2$ is uniformly bounded on the global attractor, which is known for small forces. Moreover, when the force is large and the attractor is not a fixed point, but
rather a complicated object consisting of points on complete bounded trajectories that may not be regular, the determining wavenumber $\lb_u$ from \cite{CDK} still enjoys the following
pointwise bound
\begin{equation} \label{eq:lbboundIntro}
\lb_u(t) \lesssim  \frac{\| \nabla u(t)\|_{L^2}^2}{\nu^2}.
\end{equation}
Note that this bound is optimal (from the physical point of view) in the case of extreme intermittency, i.e., when there is only one eddy at each dyadic scale. Indeed,
taking into account intermittency, Kolmogorov's dissipation wavenumber reads
\begin{equation} \label{eq:kappa_dIntro}
\kappa_\mathrm{d} := \left(\frac{\varepsilon }{\nu^3} \right)^{\frac{1}{d+1}}, \qquad \text{where} \qquad  \varepsilon := \lambda_0^{d}\nu\<\|\nabla u\|_{L^2}^2\>
=\frac{\lambda_0^{d}\nu}{T} \int_t^{t+T} \|\nabla u\|_{L^2}^2 \, d\tau.
\end{equation}
Combined with \eqref{eq:lbboundIntro}, this gives $\< \lb_u \> \lesssim \kappa_{\mathrm d}$ when $d=0$.
Here $d\in [0,3]$ is the intermittency dimension that measures the average number of eddies at various scales. Roughly speaking, the number of eddies at the lengthscale $l$
is proportional to $l^{-d}$ (see \cite{CSint} for precise mathematical definitions of active volumes, eddies, and and their relations to intermittency). In this paper we adopt an approach used in
\cite{CDK, CFr, CSr} and define the itermittency dimension $d$ through  the average level of saturation of Bernstein's inequality (see Section 3 for the precise definition).

As experimental and numerical evidence suggests, turbulent flows do not deviate much from Kolmogorov's regime where $d=3$, i.e., eddies occupy the whole region. For instance, $d \approx 2.7$ was observed in
a direct numerical simulation performed by Kaneda et al. \cite{kaneda} on the Earth Simulator.
In \cite{CDK} it was shown that one can improve \eqref{eq:lbboundIntro} for $d>0$, but such an improvement was not enough to conclude that the average determining
wavenumber was bounded by $\kappa_d$.  For instance,  in the case $d=3$, the obtained bound was $\<\lb_u\> \lesssim \kappa_{\mathrm d} ^{2+}$, which suggested that
the definition of $\lb_u$ was not optimal in the physically relevant regime. In this paper we complement the result of \cite{CDK} by 
focusing on the region $d \in [\delta,3]$, $\delta>0$,
and finding a different determining wavenumber $\lb_u$ that enjoys the optimal
bound $\<\lb_u\> \lesssim \kappa_{\mathrm{d}}$ (modulo a logarithmic correction in the case $d=3$).

We define the determining wavenumber in the following way:
\[
\lb_{u}(t):=\min\{\lambda_q:(L\lambda_{p-q})^{\frac{\delta-1}{2}}\lambda_{q}^{-1}\|u_p\|_{L^\infty}<c_0\nu ,~\forall p>q~\text{and}~ \lambda_q^{-2}\|\nabla u_{\leq q}\|_{L^\infty}<c_0\nu,~q\in \mathbb{N} \},
\]
where $0<\delta \leq3$ is a fixed (small) parameter, and $c_0$ is an adimensional constant that depends only on $\delta$. In fact, $c_0 \to 0$ as $ \delta \to 0$.
Here $\lambda_q =\frac{2^q}{L}$, $L$ is the size of the torus, $u_{\leq q}= \sum_{p=-1}^q u_q$, and $u_q = \Delta_q u$ is the Littlewood-Paley projection of $u$ (see Section~\ref{sec:pre}). Note that a convention $\min \emptyset = \infty$ is adopted in the
definition of $\lb_{u}(t)$.

Now we are ready to state our main result.

\begin{Theorem}
Let $u(t)$ and $v(t)$ be complete (ancient) bounded in $L^2$  Leray-Hopf solutions (i.e., solutions on the global attractor or pullback attractor).
Let $\lb(t):=\max\{\lb_{u}(t), \lb_{v}(t)\}$ and $Q(t)$ be such that $\lb(t)=\lambda_{Q(t)}$.
If
\begin{equation} \label{eq:dm-condition}
u(t)_{\leq Q(t)}=v(t)_{\leq Q(t)}, \qquad \forall t<0,
\end{equation}
then
\[
u(t) = v(t), \qquad \forall t \leq 0.
\]   
\end{Theorem}

The dissipation wavenumber $\Lambda_u$ enjoys the following bound:
\[
\<\lb_u\> - \l_0 \leq C_{\delta, d} \kappa_{\mathrm d}  \leq
C_{\delta, d}\kappa_0 G^{\frac{2}{d+1}}\left(\frac{1}{\nu T\kappa_0^2} + 1\right)^{\frac{1}{d+1}},
\]
for all complete bounded in $L^2$ Leray-Hopf solutions with $ d \in [\delta, 3)$.
Here $C_{\delta,d}$ is an adimensional constant that blows up when $\delta \to 0$ or $ d \to 3$.
The bound is also written in terms on the adimensional
Grashof number defined as $G:=\|f\|_{H^{-1}}/(\nu^2 \kappa_0^{1/2})$ in the autonomous case (see \eqref{eq:Gdefin} for the nonautonomous case).

In Kolmogorov's regime where $d=3$ we also obtain the optimal bound, but with a logarithmic correction:
\[
\left< \frac{\lb_u-\l_0}{(\log (\lb_u/ \l_0))^{\frac{1}{4}}} \right> \leq \widetilde C_{\delta}\kappa_{\mathrm d} \leq \widetilde C_{\delta} \kappa_0 G^{\frac{1}{2}}\left(\frac{1}{\nu T\kappa_0^2} + 1\right)^{\frac{1}{4}},
\]
for all complete bounded in $L^2$ Leray-Hopf solutions with $ d =3$.
Here $\widetilde C_{\delta}$ is an adimensional constant that depends only on the parameter $\delta$ in the definition of $\lb$. Again, $\widetilde C_{\delta} \to \infty$ as $\delta \to 0$.

\bigskip

\section{Preliminaries}
\label{sec:pre}

\subsection{Notation}
\label{sec:notation}
We denote by $A\lesssim B$ an estimate of the form $A\leq C B$ with
some absolute constant $C$, by $A\sim B$ an estimate of the form $C_1
B\leq A\leq C_2 B$ with some absolute constants $C_1$, $C_2$, and by $A\lesssim_r B$ an estimate of the form $A\leq C_r B$ with
some adimentional constant $C_r$ that depends only on the parameter $r$. We write $\|\cdot\|_p=\|\cdot\|_{L^p}$, and $(\cdot, \cdot)$ stands for the $L^2$-inner product.
We will also use $\< \cdot \>$ for time averages:
\[
\<g\>(t) := \frac{1}{T}\int_t^{t+T} g(\tau) \, d\tau, 
\]
for some fixed $T>0$.

\subsection{Littlewood-Paley decomposition}
\label{sec:LPD}
The techniques presented in this paper rely strongly on the Littlewood-Paley decomposition, which recall
here briefly. For a more detailed description on this theory we refer the readers to the books by Bahouri, Chemin and Danchin \cite{BCD} and Grafakos \cite{Gr}. 

We denote $\lambda_q=\frac{2^q}{L}$ for integers $q$. A nonnegative radial function $\chi\in C_0^\infty(\R^3)$ is chosen such that 
\begin{equation} \label{eq:xi}
\chi(\xi):=
\begin{cases}
1, \ \ \mbox { for } |\xi|\leq\frac{3}{4}\\
0, \ \ \mbox { for } |\xi|\geq 1.
\end{cases}
\end{equation}
Let 
\[
\varphi(\xi):=\chi(\xi/2)-\chi(\xi)
\]
and
\begin{equation}\notag
\varphi_q(\xi):=
\begin{cases}
\varphi(2^{-q}\xi)  \ \ \ \mbox { for } q\geq 0,\\
\chi(\xi) \ \ \ \mbox { for } q=-1,
\end{cases}
\end{equation}
so that the sequence of $\varphi_q$ forms a dyadic partition of unity. Given a tempered distribution vector field $u$ on $\T^3 =[0,L]^3$ and $q \ge -1$, an integer, the $q$th Littlewood-Paley projection of $u$ is given by 
\[
  u_q(x) := \Delta_q u(x) := \sum_{k\in\Z^3}\hat{u}(k)\phi_q(k)e^{i\frac{2\pi}{L} k \cdot x},
\]
where $\hat{u}(k)$ is the $k$th Fourier coefficient of $u$. Note that $u_{-1} = \hat{u}(0)$. Then
\[
u=\sum_{q=-1}^\infty u_q
\]
in the distributional sense. We define the $H^s$-norm in the following way:
\[
  \|u\|_{H^s} := \left(\sum_{q=-1}^\infty\lambda_q^{2s}\|u_q\|_2^2\right)^{1/2},
\]
for each $u \in H^s$ and $s\in\R$. Note that $\|u\|_{H^0} \sim \|u\|_{L^2}$. To simplify the notation, we denote
\bg\notag
u_{\leq Q}:=\sum_{q=-1}^Qu_q, \quad u_{(P,Q]}:=\sum_{q=P+1}^Qu_q, \quad \tilde{u}_q := u_{q-1} + u_q + u_{q+1}.
\ed

\subsection{Bernstein's inequality and Bony's paraproduct}
\label{sec-para}

Here we recall useful properties for the dyadic blocks of the Littlewood-Paley decomposition. The first one is the following inequality:
\begin{Lemma}\label{le:bern}(Bernstein's inequality) 
Let $n$ be the spacial dimension and $r\geq s\geq 1$. Then for all tempered distributions $u$, 
\bg\label{Bern}
\|u_q\|_{r}\lesssim \lambda_q^{n(\frac{1}{s}-\frac{1}{r})}\|u_q\|_{s}.
\ed
\end{Lemma}

Secondly, we will use the following version of Bony's paraproduct formula:
\begin{equation}\notag
\begin{split}
\Delta_q(u\cdot\nabla v)=&\sum_{|q-p|\leq 2}\Delta_q(u_{\leq{p-2}}\cdot\nabla v_p)+
\sum_{|q-p|\leq 2}\Delta_q(u_{p}\cdot\nabla v_{\leq{p-2}})\\
&+\sum_{p\geq q-2} \Delta_q(\tilde u_p \cdot\nabla v_p).
\end{split}
\end{equation}

\subsection{Weak solutions and energy inequality}
A weak solution $u(t)$ of \eqref{nse} on $[0,\infty)$ is an $L^2(\mathbb{T}^3)$ valued function in the class $u \in C([0,\infty); L^2_{\mathrm{w}}) \cap L_{\mathrm{loc}}^2(0,\infty; H^1)$ that satisfies \eqref{nse} in the sense of distributions.
A Leray-Hopf solution $u(t)$ is a weak solution satisfying the energy inequality 
\begin{equation} \label{eq:EI-sec3}
\frac{1}{2}\|u(t)\|_2^2 \leq \frac{1}{2}\|u(t_0)\|_2^2 - \nu\int_{t_0}^{t} \|\nabla u(\tau)\|_2^2\, d\tau + \int_{t_0}^{t} (f,u)\, d\tau,
\end{equation}
for almost all $t_0> 0$ and all $t >t_0$. A Leray solution $u(t)$ is a Leray-Hopf solution satisfying the above energy inequality for $t_0=0$ and all $t>t_0$. A complete
Leray-Hopf
solution $u(t)$ is an $L^2(\mathbb{T}^3)$ valued function on $(-\infty,\infty)$, such that $u(\cdot - t)|_{[0,\infty)}$ is a Leray-Hopf solution for all $t \in \mathbb{R}$.

\bigskip

\section{Global attractor, pullback attractor, and Kolmogorov's wavenumber}

In the case of a time-independent force $f$ it can be shown that the energy inequality \eqref{eq:EI-sec3} implies the  existence of an absorbing ball
\[
B:= \{ u\in L^2(\T^3): \|u\|_2 \leq R \}.
\]
Here the radius $R$ is such that
\[
R > \nu \kappa_0^{-1/2} G,
\]
where  $\kappa_0=2\pi\l_0 =2\pi/L$ and $G$ is the adimensional Grashof number
\[
G:=\frac{\|f\|_{H^{-1}}}{\nu^2 \kappa_0^{1/2}}.
\]
Note that the absorbing ball $B$ is for all the Leray solutions, i.e., the ones that satisfy the energy inequality starting from $0$.
More precisely,  for any Leray solution $u(t)$ there exists $t_0$, depending only on $\|u(0)\|_2$, such that
\[
u(t) \in B \qquad \forall t>t_0. 
\]
However, when we restrict the dynamics to
the absorbing ball, we consider Leray-Hopf solutions to define the evolutionary system and the global attractor. The Leray-Hopf solutions are weak solutions satisfying the energy inequality
starting from almost all time (but not necessarily $0$). Hence, a restriction of a Leray-Hopf solution to a smaller time interval is also a Leray-Hopf solution.
See \cite{CK} for a more detailed discussion. 

The existence of the weak global attractor $\A$ was proved in \cite{FT,FMRT}. It has the following structure:
\[
\A = \{u(0): u(\cdot) \text{ is a complete bounded  Leray-Hopf solution to the NSE}\}.
\]
The attractor $\A \subset B$ is the  $L^2$-weak omega limit of $B$, and it is the  minimal $L^2$-weakly closed weakly attracting set (see \cite{C,CF}).

In the case of a time-dependent force $f=f(t)$, a relevant object describing the long-time dynamics is a pullback attractor, whose existence was proved
in \cite{CK}. In the nonautonomous case, there exists an absorbing ball for all the Leray solutions, whose radius $R$ is such as $R > \nu \kappa_0^{-1/2} G$,
just as in the autonomous case, but the Grashof number is
\begin{equation} \label{eq:Gdefin}
G=\frac{T^{\frac12}\kappa_0^{\frac12}\|f\|_{L^2_b(T)}}{\nu^{\frac 32}(1-e^{-\nu\kappa_0^2 T})^\frac12}.
\end{equation}
Here it is assumed that $f$ is translationally bounded in $L^2_{loc}(\mathbb{R},H^{-1})$ and
\[
\|f\|^2_{L^2_b(T)}:=\sup_{t\in R} \frac{1}{T}\int_t^{t+T}\|f(\tau)\|_{H^{-1}}^2 \, d\tau.
\]
The pullback attractor is defined as the  minimal weakly closed weakly pullback attracting set for all Leray-Hopf solutions in the absorbing ball. It is the
weak pullback omega limit of $B$, and it has the following structure (see \cite{CK}):
\[
\A(t) = \{u(t): u(\cdot) \text{ is a complete bounded  Leray-Hopf solution to the NSE}\}.
\]

Let $u(t)$ be a complete bounded  Leray-Hopf solution to the NSE. Then the energy inequality \eqref{eq:EI-sec3} implies
\[
\begin{split}
0\leq\|u(t+T)\|_2^2 &\leq \limsup_{\tau \to t+}\|u(\tau)\|_2^2 - 2\nu\int_{t}^{t+T} \|\nabla u(\tau)\|_2^2\, d\tau + 2\int_{t}^{t+T} (f,u)\, d\tau\\
&\leq  \nu^2 \kappa_0^{-1} G^2 - \nu\int_{t}^{t+T} \|\nabla u(\tau)\|_2^2\, d\tau + \frac{1}{\nu}\int_{t}^{t+T} \|f\|_{H^{-1}}^2\, d\tau.\\
\end{split}
\]
Therefore
\begin{equation} \label{eq:bound G}
\begin{split}
\<\|\nabla u\|_2^2 \> := \frac{1}{T}\int_{t}^{t+T} \|\nabla u(t)\|_2^2\, dt &\leq \frac{\nu G^2}{T\kappa_0} 
+ \kappa_0\nu^2G^2.
\end{split}
\end{equation}
We can now connect this to Kolmogorov's dissipation wavenumber defined as
\begin{equation} \label{eq:kdeps-inermit}
\kappa_\mathrm{d} := \left(\frac{\varepsilon }{\nu^3} \right)^{\frac{1}{d+1}}, \qquad  \varepsilon := \nu\l_0^d\<\|\nabla u\|_2^2\>,
\end{equation}
where $d$ is the intermittency dimension and $\el$ is  average energy dissipation rate per unit active volume (i.e., the volume occupied by eddies). 
In order to define $d$, first note that
\begin{equation} \label{eq:BinInt}
\l_0^{3}\l_q^{-1} \|u_q\|_2^2 \leq  \l_q^{-1} \|u_q\|_\infty^2 \leq C_{\mathrm{B}} \l_q^{2} \|u_q\|_2^2,
\end{equation}
due to Bernstein's inequality. Here $C_B$ is an absolute constant (which depends on the choice of $\chi(\xi)$ in \eqref{eq:xi}). The intermittency dimension $d$ is defined as
\begin{equation} \label{eq:intermdef}
d:= \sup\left\{s\in \mathbb{R}: \left<\sum_{q}\l_q^{-1+s} \|u_q\|_\infty^2 \right> \leq C_{\mathrm{B}}^{3-s} \l_0^{s}\left<\sum_{q}\l_q^{2} \|u_q\|_2^2 \right> \right\},
\end{equation}
for $u \not\equiv 0$, and $d=3$ for $u \equiv 0$ on $[t, t+T]$.
Thanks to \eqref{eq:BinInt} and the fact that $\<\sum_{q}\l_q^{2} \|u_q\|_2^2 \><\infty$ , we have $d \in [0,3]$ and
\[
\left<\sum_{q}\l_q^{-1+d} \|u_q\|_\infty^2 \right> = C_{\mathrm{B}}^{3-d} \l_0^{d}\left<\sum_{q}\l_q^{2} \|u_q\|_2^2 \right>.
\]

The intermittency dimension $d$, defined in terms of  a level of saturation of Bernsten's inequality (see \cite{CSr,CSint} for similar definitions), measures the number
of eddies at various scales.
The case $d=3$ corresponds to Kolmogorov's regime where at each scale the eddies occupy the whole region. Note that $d=d(u,t)$ and
$\kappa_\mathrm{d}=\kappa_\mathrm{d}(u,t)$, defined for each individual trajectory,  are functions of time. We can also define their global analogs as
\[
D:=\inf_{u \in \mathcal{E}, t\in \mathbb{R}} d(u,t), \qquad K_{\mathrm{d}}:= \sup_{u \in \mathcal{E}, t\in \mathbb{R}} \kappa_{\mathrm{d}}(u,t).
\]
Here $\mathcal{E}$ is a family of all complete bounded  Leray-Hopf solution to the NSE.

Finally, thanks to the bound \eqref{eq:bound G},
\[
\kappa_\mathrm{d} = \left\< \frac{\l_0^d}{\nu^2 }\|\nabla u\|_2^2 \right\>^{\frac{1}{d+1}} \leq
(2 \pi)^{-\frac{d}{d+1}}\kappa_0G^{\frac{2}{d+1}}\left(\frac{1}{\nu T\kappa_0^2} + 1\right)^{\frac{1}{d+1}}. 
\]
Also, taking the supremum over all $u \in \mathcal{E}$ and $t\in \mathbb{R}$, we obtain 
\[
K_{\mathrm{d}} \leq \kappa_0G^{\frac{2}{D+1}}\left(\frac{1}{\nu T\kappa_0^2} + 1\right)^{\frac{1}{D+1}},
\]
provided $G\geq 1$.

\section{Proof of the main result}
\label{sec:pf}

Let $u(t)$ and $v(t)$ be complete bounded in $L^2$  Leray-Hopf solutions. Denote $w:=u-v$, which satisfies the equation
\begin{equation} \label{eq-w}
w_t+u\cdot\nabla w+w\cdot\nabla v=-\nabla p'+\nu \Delta w
\end{equation}
in the sense of  distributions. Here $p'$ stands for the difference of the pressures.

Recall the definition of the determining wavenumber:
\[
\lb_{u}(t)=\min\{\lambda_q:(L\lambda_{p-q})^{\sigma}\lambda_{q}^{-1}\|u_p\|_{L^\infty}<c_0\nu ,~\forall p>q~\text{and}~ \lambda_q^{-2}\|\nabla u_{\leq q}\|_{L^\infty}<c_0\nu,~q\in \mathbb{N} \},
\]
where $\sigma = (\delta-1)/2$. 
Let $\lb(t):=\max\{\lb_{u}(t), \lb_{v}(t)\}$ and $Q(t)$ be such that $\lb(t)=\lambda_{Q(t)}$. By our assumption, $w_{\leq Q(t)}(t)\equiv 0$.
Recall that $0<\delta \leq 3$, i.e., $-1/2<\sigma\leq 1$. Let
\[
s =\min\left\{\textstyle -\frac12+\frac{\delta}{4},0\right\}.
\]
Then straightforward computations give $-1-\sigma <s <\sigma\leq 1$.

Multiplying (\ref{eq-w}) by $\ls\Delta^2_qw$, integrating (i.e., using $\ls\Delta^2_qw$ as a test function in the weak formulation), and adding up for all $q\geq -1$ yields 
\begin{equation}\label{w2}
\begin{split}
\frac{1}{2}\|w(t)\|_{H^s}^2- \frac{1}{2}\|w(t_0)\|_{H^s}^2 +\nu &\int_{t_0}^t\|w\|_{H^{1+s}}^2 \, d\tau\\
\leq & \int_{t_0}^t \sum_{q\geq -1}\ls\left|\int_{\T^3}\Delta_q(w\cdot\nabla v) w_q\, dx\right|\, d\tau\\
&+\int_{t_0}^t\sum_{q\geq -1}\ls\left|\int_{\T^3}\Delta_q(u\cdot\nabla w) w_q\, dx\right|\, d\tau,\\
=:&\int_{t_0}^t I\, d\tau+ \int_{t_0}^t J\, d\tau.
\end{split}
\end{equation}

We first decompose $I$ using Bony's paraproduct as mentioned in Subsection \ref{sec-para},
\begin{equation}\notag
\begin{split}
I\leq& \sum_{q\geq -1}\ls\sum_{|q-p|\leq 2}\left|\int_{\T^3}\Delta_q(w_{\leq{p-2}}\cdot\nabla v_p) w_q\, dx\right|\\
&+\sum_{q\geq -1}\ls\sum_{|q-p|\leq 2}\left|\int_{\T^3}\Delta_q(w_{p}\cdot\nabla v_{\leq{p-2}})w_q\, dx\right|\\
&+\sum_{q\geq -1}\ls\sum_{p\geq q-2} \left|\int_{\T^3}\Delta_q(\tilde w_p\cdot\nabla v_p)w_q\, dx\right|\\
=:&I_{1}+I_{2}+I_{3}.
\end{split}
\end{equation}
It follows from H\"older's inequality that
\[
\begin{split}
  I_{1} 
& \leq \sum_{q>Q}\sum_{\substack{|q-p|\leq 2\\p>Q+2}}\ls\int_{\T^3}|\Delta_q(w_{\leq{p-2}}\cdot\nabla v_p) w_q|\, dx\\
 &\lesssim  \sum_{q>Q}\sum_{\substack{|q-p|\leq 2\\p>Q+2}}\ls\|w_{(Q, p-2]}\|_{2}\lambda_p\|v_p\|_\infty\|w_q\|_2.
\end{split}
\]
Using the definition of $\lb$, Young's inequality and Jensen's inequality, we obtain
\[
\begin{split}
 I_{1}  &\lesssim  c_0\nu\sum_{q>Q}\sum_{\substack{|q-p|\leq 2\\p>Q+2}}\ls\lb^{1+\sigma}\lambda_p^{1-\sigma}\|w_q\|_2\sum_{Q<p'\leq p-2}\|w_{p'}\|_{2} \\
 &\lesssim  c_0\nu\sum_{q>Q}\lambda_q^{1+s}\|w_q\|_2\left(\sum_{Q<p'\leq q}\lambda_{p'}^{1+s}\|w_{p'}\|_2\lambda_{p'}^{-1-s}\lambda_{q}^{s-\sigma}\lambda_Q^{1+\sigma}\right)\\
&\lesssim  c_0\nu \sum_{q>Q}\lambda_q^{1+s}\|w_q\|_2\left(\sum_{Q<p'\leq q}\lambda_{p'}^{1+s}\|w_{p'}\|_2(L\l_{q-p'})^{s-\sigma}\right),
\end{split}
\]
where we used $\sigma\geq -1$ and $s<\sigma$.
Now using Young's inequality and Jensen's inequality, we conclude
\[
\begin{split}
I_1&\lesssim  c_0\nu \sum_{q>Q}\lambda_q^{2+2s}\|w_q\|_2^2+c_0\nu \sum_{q>Q}\left(\sum_{Q<p'\leq q}\lambda_{p'}^{1+s}\|w_{p'}\|_2(L\l_{q-p'})^{s-\sigma}\right)^2\\
&\lesssim  c_0\nu \sum_{q>Q}\lambda_q^{2+2s}\|w_q\|_2^2+c_0\nu \sum_{q>Q}\sum_{Q<p'\leq q}\lambda_{p'}^{2+2s}\|w_{p'}\|_2^2(L\l_{q-p'})^{s-\sigma}\\
&\lesssim  c_0\nu \sum_{q>Q}\lambda_q^{2+2s}\|w_q\|_2^2+c_0\nu \sum_{p'>Q}\lambda_{p'}^{2+2s}\|w_{p'}\|_2^2\sum_{q\geq p'}(L\l_{q-p'})^{s-\sigma}\\
&\lesssim  c_0\nu \|\nabla^{1+s} w\|_2^2,
\end{split}
\]
where we needed $s<\sigma$. Note that we omit adimensional constants that depend on $\delta$
throughout this proof. The precise bound on $I_1$ is
\[
I_{1} \lesssim c_0\nu \|\nabla^{1+s} w\|_2^2\left(1+(1-2^{s-\sigma})^{-1}\right).
\]
Note that $(1-2^{s-\sigma})^{-1} \to \infty$ as $\delta \to 0+$ by definitions of $\sigma$ and $s$. Because of
this we will have $c_0 \to 0$ as $\delta\to 0+$ once we choose $c_0$ at the end of the proof.
This explains why we have to avoid the case of extreme intermittency, which is covered in the companion paper \cite{CDK}.

Following a similar strategy, we have
\[
\begin{split}
  I_{2} &\leq \sum_{q>Q}\sum_{\substack{|q-p|\leq 2\\p>Q}}\ls\int_{\T^3}|\Delta_q(w_p\cdot\nabla v_{\leq p-2}) w_q| \, dx\\
&\leq  \sum_{q>Q}\sum_{\substack{|q-p|\leq 2\\p>Q}}\ls\|w_p\|_2\|\nabla v_{(Q,p-2]}\|_\infty\|w_q\|_2\\
&+\sum_{q>Q}\sum_{\substack{|q-p|\leq 2\\p>Q}}\ls\|w_p\|_2\|\nabla v_{\leq Q}\|_\infty\|w_q\|_2\\
&\equiv I_{21}+I_{22},
\end{split}
\]
where we adopt the convention that $(Q,p-2]$ is empty if $p-2\leq Q$.
Thus, the first part of the definition of $\lb$ implies
\[
\begin{split}
I_{21}
           &\lesssim \sum_{p> Q}\sum_{|q-p|\leq 2}\ls\| w_p\|_2\|w_q\|_{2}\sum_{Q<p'\leq p-2}\|\nabla v_{p'}\|_\infty\\
           &\lesssim \sum_{q> Q}\ls\| w_q\|_2^2\sum_{Q<p'\leq q+2}\lambda_{p'}\|v_{p'}\|_\infty\\
           &\lesssim c_0\nu\sum_{q> Q}\ls\| w_q\|_2^2\sum_{Q<p'\leq q+2}\lambda_{p'}^{1-\sigma}\lb^{1+\sigma}\\
           &\lesssim c_0\nu\sum_{q> Q}\l_q^{2+2s}\| w_q\|_2^2\sum_{Q<p'\leq q+2}\lambda_{p'}^{1-\sigma}\lb^{1+\sigma}\l_q^{-2}\\
           &\lesssim  c_0\nu \sum_{q> Q}\l_q^{2+2s}\| w_q\|_2^2
\end{split}
\]
where we need $\sigma\geq-1$. While the second part of the definition of $\lb$ gives
\[
\begin{split}
I_{22}
&\lesssim  \sum_{q> Q}\sum_{\substack{|q-p|\leq 2\\p>Q}}\ls \|w_p\|_2\|w_q\|_2\|\nabla v_{\leq Q}\|_\infty\\
&\lesssim  c_0\nu\sum_{q> Q}\lb_v^2\ls \|w_q\|_2^2\\
&\lesssim  c_0\nu \sum_{q> Q}\l_q^{2+2s} \|w_q\|_2^2.
\end{split}
\]

We will now  estimate $I_{3}$.  It follows from integration by parts that 
\[
\begin{split}
  I_{3}&=\sum_{q\geq -1}\ls\sum_{|q-p|\leq 2}\left|\int_{\T^3}\Delta_q(w_{p}\cdot\nabla v_{\leq{p-2}})w_q\, dx\right|\\ 
  &\leq  \sum_{q>Q}\sum_{p\geq q-2}\ls\int_{\T^3}|\Delta_q(\tilde w_p \otimes v_{p}) \nabla w_q| \, dx\\
&\leq  \sum_{p> Q}\sum_{Q<q\leq p+2}\l_q^{1+2s}\|\tilde w_p\|_2\|w_q\|_2\|v_p\|_\infty.
\end{split}
\]
By H\"older's inequality and definition of $\lb$ we have
\[
\begin{split}
I_{3}
&\lesssim  \sum_{p> Q}\|\tilde w_p\|_2\|v_{p}\|_\infty\sum_{Q<q\leq p+2}\lambda_q^{1+2s}\|w_q\|_2\\
&\lesssim  c_0\nu\sum_{p> Q}\lb^{1+\sigma}\l_p^{-\sigma}\|\tilde w_p\|_2\sum_{Q<q\leq p+2}\lambda_q^{1+2s }\|w_q\|_2\\
&\lesssim  c_0\nu\sum_{p> Q}\l_p^{1+s}\|\tilde w_p\|_2\sum_{Q<q\leq p+2}\lambda_q^{1+s}\|w_q\|_2\l_Q^{1+\sigma}\l_p^{-1-s-\sigma}\l_q^s\\
\end{split}
\]
Now we use Young's and Jensen's inequalities to infer 
\[
\begin{split}
I_{3}&\lesssim  c_0\nu \sum_{p> Q}\l_p^{1+s}\|\tilde w_p\|_2\sum_{Q<q\leq p+2}\lambda_q^{1+s}\|w_q\|_2(L\l_{q-p})^{1+s+\sigma}\\
&\lesssim  c_0\nu \sum_{p> Q}\l_p^{2+2s}\|w_p\|_2^2+c_0\nu \sum_{p>Q}\left(\sum_{Q<q\leq p+2}\lambda_q^{1+s}\|w_q\|_2(L\l_{q-p})^{1+s+\sigma}\right)^2\\
&\lesssim  c_0\nu \sum_{p> Q}\l_p^{2+2s}\|w_p\|_2^2,
\end{split}
\]
where we used $\sigma\geq-1$ and $s>-1-\sigma$.

Therefore, we have for $\sigma\geq -1$ and $-1-\sigma<s< \sigma$,
\begin{equation}\label{est-i1}
I\lesssim c_0 \nu \|\nabla^{1+s} w\|_2^2.
\end{equation}

Now applying Bony's paraproduct formula to $J$ yields
\begin{equation}\notag
\begin{split}
J=&\int_{t_0}^t \sum_{q\geq -1}\ls\left|\int_{\T^3}\Delta_q(w\cdot\nabla v) w_q\, dx\right|\, d\tau\\
\leq&\sum_{q\geq -1}\sum_{|q-p|\leq 2}\ls\left|\int_{\T^3}\Delta_q(u_{\leq{p-2}}\cdot\nabla w_p)w_q\, dx\right|\\
&+\sum_{q\geq -1}\sum_{|q-p|\leq 2}\ls\left|\int_{\T^3}\Delta_q(u_{p}\cdot\nabla w_{\leq{p-2}})w_q\, dx\right|\\
&+\sum_{q\geq -1}\sum_{p\geq q-2}\sum_{|p-p'|\leq 1}\ls\left|\int_{\T^3}\Delta_q(u_p\cdot\nabla w_{p'})w_q\, dx\right|\\
=:& J_{1}+J_{2}+J_{3}.
\end{split}
\end{equation}
We further decompose $J_{1}$ by using a commutator form 
\begin{equation}\notag
\begin{split}
J_{1}\leq& \sum_{q\geq -1}\sum_{|q-p|\leq 2}\ls\left|\int_{\R^3}[\Delta_q, u_{\leq{p-2}}\cdot\nabla] w_pw_q\, dx\right|\\
&+\sum_{q\geq -1}\ls\left|\int_{\R^3} u_{\leq q-2}\cdot\nabla w_q w_q\, dx\right|\\
&+\sum_{q\geq -1}\sum_{|q-p|\leq 2}\ls\left|\int_{\R^3}( u_{\leq{p-2}}- u_{\leq q-2})\cdot\nabla\Delta_qw_p w_q\, dx\right| \\
=&J_{11}+J_{12}+J_{13}.
\end{split}
\end{equation}
To obtain the second term  we used $\sum_{|p-q|\leq 2}\Delta_qw_p=w_q$.  In fact, we have $J_{12}=0$ since $\div\, u_{\leq q-2}=0$. In the first term, the commutator is defined as

\[[\Delta_q, u_{\leq{p-2}}\cdot\nabla] w_p:=\Delta_q(u_{\leq{p-2}}\cdot\nabla w_p)-u_{\leq{p-2}}\cdot\nabla \Delta_qw_p.\]
 It is easy to see  (see \cite{CD} for more details) that for any $1\leq r\leq \infty$,
\begin{equation}\notag
\|[\Delta_q, u_{\leq{p-2}}\cdot\nabla] w_p\|_r\\
\lesssim \|\nabla  u_{\leq p-2}\|_\infty\|w_p\|_r.
\end{equation}
Then $J_{11}$ is estimated as
\[
\begin{split}
  J_{11} &\leq  \sum_{q>Q}\sum_{\substack{|q-p|\leq 2\\p>Q}}\ls\int_{\T^3}|[\Delta_q, u_{\leq{p-2}}\cdot\nabla] w_pw_q|\, dx\\
  &\leq  \sum_{q>Q}\sum_{\substack{|q-p|\leq 2\\p>Q}}\ls\|\nabla u_{(Q,p-2]}\|_\infty\|w_p\|_2\|w_q\|_2\\
  &+\sum_{q>Q}\sum_{\substack{|q-p|\leq 2\\p>Q}}\ls\|\nabla u_{\leq Q}\|_\infty\|w_p\|_2\|w_q\|_2\\
  &\equiv J_{111}+J_{112}.
\end{split}
\]
Here
\[
\begin{split}
J_{111}
&\lesssim  \sum_{q>Q}\l_q^{2s}\|w_q\|_2^2\sum_{Q<p'\leq q}\l_{p'}\|u_{p'}\|_\infty\\
&\lesssim  c_0\nu\sum_{q>Q}\l_q^{2s}\|w_q\|_2^2\sum_{Q<p'\leq q}\lb^{1+\sigma}\l_{p'}^{1-\sigma}\\
&\lesssim  c_0\nu\sum_{q>Q}\l_q^{2+2s}\|w_q\|_2^2\sum_{Q<p'\leq q}\lb^{1+\sigma}\l_{p'}^{1-\sigma}\l_q^{-2}\\
&\lesssim  c_0\nu \sum_{q>Q}\l_q^{2+2s}\|w_q\|_2^2,
\end{split}
\]
where we used $\sigma\geq-1$. As for the second term,
using the fact that $\|\nabla u_{\leq q}\|_\infty\leq c_0\nu\lb^2$ for $q\leq Q$, we obtain 
\[
J_{112}
\lesssim  c_0\nu\lb^2\sum_{q>Q}\l_q^{2s}\|w_q\|_2^2
\lesssim  c_0\nu \sum_{q>Q}\l_q^{2+2s}\|w_q\|_2^2.
\]

The term $J_{13}$ is estimated as
\[
\begin{split}
J_{13}
&\leq \sum_{q>Q}\sum_{\substack{|q-p|\leq 2\\p>Q}}\ls\int_{\R^3}|( u_{\leq{p-2}}- u_{\leq q-2})\cdot\nabla\Delta_qw_p w_q|\, dx\\
&\lesssim \sum_{q>Q}\l_q^{1+2s}\|u_{(q-4, q]}\|_\infty\|w_q\|_2^2\\
&\lesssim \sum_{q>Q}\l_q^{1+2s}\|u_{(q-4, Q]}\|_\infty\|w_q\|_2^2+\sum_{q>Q}\sum_{\substack{q-4<p'\leq q\\p'>Q}}\l_q^{1+2s}\|u_{p'}\|_\infty\|w_q\|_2^2\\
&\equiv J_{131}+J_{132}.
\end{split}
\]
As before, we adopt the convention that $(q-4,Q]$ is empty if $q-4\geq Q$. We have 
\[
J_{131}\lesssim c_0\nu \lb\sum_{q>Q}\l_q^{1+2s}\|w_q\|_2^2\lesssim c_0\nu \sum_{q>Q}\l_q^{2+2s}\|w_q\|_2^2,
\]
and
\[
\begin{split}
J_{132}& =\sum_{q>Q}\sum_{\substack{q-4\leq p'\leq q\\p'>Q}}\l_q^{1+2s}\|u_{p'}\|_\infty\|w_q\|_2^2\\
& \lesssim c_0\nu\sum_{q>Q}\sum_{\substack{q-4\leq p'\leq q\\p'>Q}}\l_q^{1+2s}\lb^{1+\sigma}\l_{p'}^{-\sigma}\|w_q\|_2^2\\
& \lesssim c_0\nu\sum_{q>Q}\l_q^{2+2s}\|w_q\|_2^2(L\l_{Q-q})^{1+\sigma}\\
& \lesssim c_0\nu\sum_{q>Q}\l_q^{2+2s}\|w_q\|_2^2,
\end{split}
\]
where we used $\sigma\geq -1$.

Now we continue with  $J_2$:
\[
\begin{split}
  J_{2} &= \sum_{q>Q}\sum_{\substack{|q-p|\leq 2\\p>Q+2}}\ls\left|\int_{\T^3}\Delta_q(u_{p}\cdot\nabla w_{\leq{p-2}})w_q\, dx\right|\\
  &\leq \sum_{q>Q}\sum_{\substack{|q-p|\leq 2\\p>Q+2}}\ls\|u_p\|_\infty\|\nabla w_{(Q, p-2]}\|_2\|w_q\|_2.
   \end{split}
   \]
Using definition of $\lb$, Young's, and Jensen's inequalities we obtain   
 \[
   \begin{split}
  J_{2}
&\lesssim  c_0\nu\sum_{q>Q}\sum_{\substack{|q-p|\leq 2\\p>Q+2}}\ls\lb^{1+\sigma}\l_p^{-\sigma}\|w_q\|_2\|\nabla w_{(Q, p-2]}\|_2\\
&\lesssim  c_0\nu\sum_{q>Q}\lb^{1+\sigma}\l_q^{2s-\sigma}\|w_q\|_2\|\nabla w_{(Q, q]}\|_2\\
&\lesssim  c_0\nu\sum_{q>Q}\lb^{1+\sigma}\l_q^{2s-\sigma}\|w_q\|_2\sum_{Q<p'\leq q}\l_{p'}\| w_{p'}\|_2\\       
&\lesssim  c_0\nu\sum_{q>Q}\l_q^{1+s}\|w_q\|_2\sum_{Q<p'\leq q}\l_{p'}^{1+s}\| w_{p'}\|_2\l_q^{s-\sigma-1}\l_{p'}^{-s}\lb^{1+\sigma}\\     
&\lesssim  c_0\nu\sum_{q>Q}\l_q^{1+s}\|w_q\|_2\left(\sum_{Q<p'\leq q}\l_{p'}^{1+s}\| w_{p'}\|_2(L\l_{q-p'})^{s-\sigma-1}\right)\\
&\lesssim  c_0\nu\sum_{q>Q}\l_q^{2+2s}\|w_q\|_2^2+c_0\nu\sum_{q>Q}\left(\sum_{Q<p'\leq q}\l_{p'}^{1+s}\| w_{p'}\|_2(L\l_{q-p'})^{s-\sigma-1}\right)^2\\
&\lesssim  c_0\nu\sum_{q>Q}\l_q^{2+2s}\|w_q\|_2^2,
\end{split}
\]
where we used $s<\sigma+1$ and $\sigma\geq -1$.

Notice that the last term $J_{3}$ can be estimated in the same way as $I_{3}$.
Therefore we have for $\sigma\geq -1$ and $-1-\sigma<s< 1+\sigma$,
\begin{equation}\label{est-i2}
J \lesssim c_0 \nu \|\nabla^{1+s} w\|_2^2.
\end{equation}

Combining \eqref{est-i1} and \eqref{est-i2}, we conclude that for any $\delta>0$, there exists an adimensional constant $C>0$ (that depends only on $\delta$) such that
\[
I+J \leq C c_0 \nu \|\nabla^{1+s} w\|_2^2,
\]
where $s =\min\left\{\textstyle -\frac12+\frac{\delta}{4},0\right\} \leq 0$.
Choosing $c_0:=1/(2C)$ we infer from \eqref{w2} that for all $t_0\leq t$,
\[
\begin{split}
\|w(t)\|_{H^s}^2 \leq & \|w(t_0)\|_{H^s}^2 - \nu \int_{t_0}^t \|\nabla^{1+s} w\|_2^2 \, d\tau\\
\leq & \|w(t_0)\|_{H^s}^2 - \nu \k_0^{2+2s} \int_{t_0}^t \|w\|_2^2 \, d\tau,
\end{split}
\]
with $\k_0=\frac{2\pi}{L}$.
Thus 
\[
\|w(t)\|_{H^s}^2 \leq \|w(t_0)\|_{H^s}^2e^{-\nu\k_0^{2+2s}(t-t_0)}, \qquad t_0\leq t.
\]
Recall that $s \leq 0$ and hence $\|w(t)\|_{H^s} \lesssim \lambda_0^{s} \|w(t)\|_2$, which is bounded on
$\mathbb{R}$ as $w(t)$ is the difference of two complete bounded trajectories. Taking the limit as $t_0 \to -\infty$ completes the proof.

\cbdu

\section{Average determining wavenumber and  Kolmogorov's dissipation wavenumber}
\label{sec:Kolmogorov}

The goal of this section is to derive a uniform upper bound on the average determining wavenumber in the absorbing ball.
First, recall that $\lb_u(t)$ is defined as 
\[
\lb_{u}(t):=\min\{\lambda_q:(L\lambda_{p-q})^{\sigma}\lambda_{q}^{-1}\|u_p\|_{\infty}<c_0\nu ,~\forall p>q~\text{and}~ \lambda_q^{-2}\|\nabla u_{\leq q}\|_{\infty}<c_0\nu,~q\in \mathbb{N} \},
\]
where $\sigma = (\delta-1)/2$ and $c_0$ is an adimensional constant that depends only on $\delta$. Recall that $\sigma \in (-1/2, 1]$. We will drop the subscript $u$ in $\lb_u$ and define $Q $ so that $\lambda_Q = \lb$.
\begin{Lemma} \label{L:Lambda-main-estimate}
 If $\l_0 \leq \lb<\infty$, then
\begin{equation} \label{eq:Lambda-main-estimate}
(c_0 \nu)^2\lb^4 \lesssim  \|\nabla u_{\leq Q-1}\|_\infty^2 + \sup_{p\geq Q}  (L\l_{p-Q})^{2\sigma}\lb^2 \|u_p\|_\infty^2.
\end{equation}
If $\lb=\infty$, then
\[
\sup_q \l_q^\sigma \|u_q\|_\infty = \infty.
\]
\end{Lemma}
\begin{proof}
First, consider the case $\lb=\infty$. Then for every $q\in \mathbb{N}$ either
\begin{equation} \label{1st-cond-in-L}
\sup_{p>q} (L\l_{p-q})^\sigma \l_{q}^{-1}\|u_p\|_{\infty}\geq c_0\nu ,
\end{equation}
or
\begin{equation} \label{2nd-cond-in-L}
\lambda_q^{-2}\|\nabla u_{\leq q}\|_{\infty}\geq c_0\nu.
\end{equation}

If \eqref{1st-cond-in-L} is satisfied for infinitely many $q\in \mathbb{N}$, then
\[
\limsup_{q\to \infty} \sup_{p>q} \l_q^{-\sigma-1} (L\l_{p})^{\sigma}\|u_{p}\|_{\infty}\geq c_0\nu.
\]
Since $\sigma > -1$, this immediately 
implies that $\sup_q \l_q^\sigma \|u_q\|_\infty = \infty$.

If \eqref{2nd-cond-in-L} is satisfied for infinitely
many $q\in \mathbb{N}$, then
\[
\limsup_{q\to \infty} \lambda_q^{-2}\|\nabla u_{\leq q}\|_{\infty}\geq c_0\nu.
\]
On the other hand, since $\sigma \leq 1$,
\[
\begin{split}
\lambda_q^{-2}\|\nabla u_{\leq q}\|_{\infty} &\lesssim \l_q^{-2}\sum_{p\leq q} \l_p\|u_p\|_\infty\\
&=\l_q^{-\sigma-1}\sum_{p\leq q} (L\lambda_{q-p})^{\sigma-1}\l_p^\sigma\|u_p\|_\infty\\
&\leq \l_q^{-\sigma-1} \sup_{p\leq q} \l_p^\sigma\|u_p\|_\infty.
\end{split}
\]
Hence, since $-\sigma-1 <0$,  $\sup_q \l_q^\sigma \|u_q\|_\infty = \infty$.

Now if $\l_0<\lb(t) <\infty$,
then both conditions in the definition of $\lb$ are satisfied for $q=Q$, but one of the conditions is not satisfied for $q = Q-1$, i.e.,
\begin{equation}\label{eq:alt1}
2^{(p-Q+1)\sigma}\l_{Q-1}^{-1}\|u_p\|_\infty \geq c_0\nu, \qquad \text{for some} \qquad p\geq Q,
\end{equation}
or
\begin{equation} \label{eq:alt2}
 \|\nabla u_{\leq Q-1}\|_\infty \geq c_0 \nu\lambda_{Q-1}^2={\textstyle \frac{1}{4}}c_0 \nu\lb^2.
\end{equation}
Thus we have
\[
(c_0\nu)^2 \lb^4 \leq  16 (\l_{p-Q}L)^{2\sigma}\lb^2 \|u_p\|_\infty^2, \qquad \text{for some} \qquad p\geq Q, 
\]
or
\[
(c_0\nu)^2 \lb^4 \leq 16 \|\nabla u_{\leq Q-1}\|_\infty^2.
\]
Hence, adding the right hand sides, we obtain
\eqref{eq:Lambda-main-estimate}.
\end{proof}
We will now consider the average determining wavenumber
\[
\<\lb\> := \frac{1}{T}\int_t^{t+T} \lb(\tau) \, d\tau,
\]
and compare it to Kolmogorov's dissipation wavenumber defined as
\begin{equation} \label{eq:kdeps-inermit}
\kappa_\mathrm{d} := \left(\frac{\varepsilon }{\nu^3} \right)^{\frac{1}{d+1}}, \qquad  \varepsilon := \nu\l_0^d\<\|\nabla u\|_2^2\>
= \frac{\nu\l_0^d}{T}\int_t^{t+T} \|\nabla u(\tau)\|_2^2 \, d\tau,
\end{equation}
where $d\in[0,3]$ is the intermittency dimension and $\el$ is  average energy dissipation rate per unit active volume (i.e., the volume occupied by eddies).
Recall from the definition of intermittency \eqref{eq:intermdef} that
\begin{equation} \label{eq:intermdef}
\<\sum_{q\leq Q}\l_q^{-1+d} \|u_q\|_\infty^2 \> \lesssim \l_0^{d}\<\sum_{q\leq Q}\l_q^{2} \|u_q\|_2^2 \>.
\end{equation}
The case $d=3$ corresponds to Kolmogorov's regime where at each scale the eddies occupy the whole region, and $d=0$ is the case of extreme intermittency.

Consider now a solution $u$ for which $d \geq \delta$, i.e., $d \geq 2\sigma +1$.
Then whenever $\lb_u(t)$ is finite, we can use \eqref{eq:Lambda-main-estimate}
in Lemma~\ref{L:Lambda-main-estimate} and Jensen's
inequality to get
\[
\begin{split}
(\lb - \l_0)^{d+1}
& \lesssim  \frac{\lb^{d-3}}{(c_0 \nu)^2}\left( \|\nabla u_{\leq Q-1}\|_\infty^2 + \sup_{q\geq Q}  (L\l_{q-Q})^{2\sigma} \lb^2\|u_q\|_\infty^2\right)\\
&\lesssim  \frac{1}{\nu^2}\left(\sum_{q\leq Q-1} \l_q^{(d-1)/2}\|u_q\|_\infty (L\lambda_{Q-q})^{(d-3)/2} \right)^2
+ \frac{\lb^{d-1}}{\nu^2}\sup_{q\geq Q}  (L\l_{q-Q})^{2\sigma} \|u_q\|_\infty^2\\
&\lesssim_d \frac{1}{\nu^2} \sum_{q\leq Q-1} \l_q^{d-1}\|u_q\|_\infty^2 + \frac{1}{\nu^2}\sup_{q\geq Q}  (L\l_{q-Q})^{2\sigma-d+1} \l_q^{d-1}\|u_q\|_\infty^2\\
&\lesssim_d \frac{1}{\nu^2} \sum_{q} \l_q^{d-1}\|u_q\|_\infty^2.
\end{split}
\]
If $\lb=\infty$, this inequality is also true. Indeed, in this case Lemma~\ref{L:Lambda-main-estimate} implies
\[
\sum_{q} \l_q^{d-1}\|u_q\|_\infty^2 \geq \sum_{q} \l_q^{2\sigma}\|u_q\|_\infty^2 =\infty.
\]
Then thanks to Jensen's inequality,
\[
\begin{split}
\<\lb\>-\l_0
& \lesssim \<(\lb - \l_0)^{d+1}\>^{\frac{1}{d+1}}\\
&\lesssim_d \left\<\frac{1}{\nu^2} \sum_{q} \l_q^{d-1}\|u_q\|_\infty^2 
\right\>^{\frac{1}{d+1}}.
\end{split}
\]
Now using \eqref{eq:intermdef} we conclude that
\[
\begin{split}
\<\lb\>-\l_0
&\lesssim_d \left\<\frac{1}{\nu^2} \sum_{q\leq Q} \l_q^{d-1}\|u_q\|_\infty^2 \right\>^{\frac{1}{d+1}}\\
&\lesssim \left\< \frac{\l_0^d}{\nu^2}\sum_{q\leq Q} \l_q^{2}\|u_q\|_2^2 \right\>^{\frac{1}{d+1}}\\
&\lesssim  \left\< \frac{\nu \l_0^d}{\nu^3 }\|\nabla u\|_2^2 \right\>^{\frac{1}{d+1}}\\
&= \kappa_\mathrm{d}
\end{split}
\]

Consider now  Kolmogorov's regime where $d=3$. Then a similar computation yields
\[
\begin{split}
\left< \frac{\lb-\l_0}{(\log (\lb/ \l_0))^{\frac{1}{4}}} \right>
& \lesssim \left<\frac{(\lb - \l_0)^{4}}{Q}\right>^{\frac{1}{4}}\\
& \lesssim \left< Q\left(\frac{1}{c_0 \nu Q} \sum_{q\leq Q} \|\nabla u_q\|_\infty\right)^2 + \sup_{p\geq Q}  (L\l_{p-Q})^{2\sigma}\lb^2 \|u_p\|_\infty^2 \right>^{\frac{1}{4}}\\
&\lesssim \left\< \frac{1}{\nu^2}\sum_{q\leq Q} \l_q^2 \|u_q\|_\infty^2  + \frac{1}{\nu^2}\sup_{q\geq Q}  (L\l_{q-Q})^{2\sigma-2} \l_q^{2}\|u_q\|_\infty^2
 \right\>^{\frac{1}{4}}\\
&\lesssim \left\< \frac{\l_0^3}{\nu^2}\sum_{q} \l_q^{2}\|u_q\|_2^2 \right\>^{\frac{1}{4}}\\
&\lesssim \kappa_\mathrm{d}.
\end{split}
\]

\section*{Acknowledgments} 
The authors thank the anonymous referee for careful reading the manuscript and
constructive comments.

\end{document}